\documentclass[a4paper,12pt,english]{article}
\usepackage[utf8x]{inputenc}
\usepackage{amsmath, amsthm, amssymb, graphicx, makeidx, mathrsfs, enumerate}
\usepackage{setspace}
\usepackage{hyperref}
\usepackage{multirow}
\usepackage[T1]{fontenc}
\usepackage{babel}
\usepackage{authblk}

\newcommand{\Z}{\mbox{$\mathbb Z$}}	
\newcommand{\Q}{\mbox{$\mathbb Q$}}	

\usepackage{blindtext}
\usepackage{sectsty}
\sectionfont{\centering}

\usepackage{vmargin}
\setmarginsrb     { 1.0in}  
                        { 1.0in}  
                        { 1.0in}  
                        { 1.0in}  
                        {  20pt}  
                        {0.25in}  
                        {   9pt}  
                        { 0.3in}  
\newtheorem{theorem}{Theorem}[section]

\newtheorem{corollary}[theorem]{Corollary}
\theoremstyle{definition}
\newtheorem{definition}[theorem]{Definition}
\newtheorem{notation}[theorem]{Notation}
\newtheorem{example}[theorem]{Example}

\theoremstyle{remark}
\newtheorem{remark}[theorem]{Remark}

\numberwithin{equation}{section}

\newcommand{\x}{\mbox{$ \mathfrak{d} $}}
\newcommand{\y}{\ell}
\makeindex

\usepackage{etoolbox}
\patchcmd{\thebibliography}
  {\settowidth}
  {\setlength{\parsep}{0pt}\setlength{\itemsep}{0pt plus 0.1pt}\settowidth}
  {}{}

\begin{document}
\addcontentsline{toc}{chapter}{\bf Notation}
\addcontentsline{toc}{chapter}{\bf Asymptotic Notation for Runtime Analysis}
\pagenumbering{arabic}

\date{}
\renewcommand\Authands{}

\title{\large{\bf{\textsc{On the irreducible factors of a polynomial II}}}}
\author{Anuj Jakhar,  ~Srinivas Kotyada}
\affil[]{Institute of Mathematical Sciences, HBNI, CIT Campus, Taramani, Chennai - 600113, Tamil Nadu, India.\\anujjakhar@imsc.res.in, srini@imsc.res.in}

\maketitle
\begin{center}
{\large{\bf {\textsc{Abstract}}}}
\end{center}

We give a lower bound for the degree of an irreducible factor of a given polynomial. This improves and generalizes the results obtained in \cite[On the irreducible factors of a polynomial, to appear in Proc. Amer. Math. Soc., 2019]{Jak}.

\bigskip


\bigskip

\section{\textsc{Introduction}}
In \cite[Theorem 1.1]{Jak}, the following result was proved for polynomials having integer coefficients.
\begin{theorem}\label{a1}
Let $p$ be a prime number and let $f(x) = a_nx^n + a_{n-1}x^{n-1} + \cdots + a_0$  be a polynomial with integer coefficients. Suppose that $p$ does not divide $a_s$ for some $s\leq n$, and that $a_j \neq 0$ for some $j$ with $0\leq j<s$. For $0\leq i <s$, let  $r_i$ be the largest positive integer such that $p^{r_i}$ divides $a_i$ (where $r_i = \infty$ if $a_i = 0$). Let $k (<s)$ be the smallest non-negative integer such that $\min\limits_{0\leq i<s}\frac{r_i}{s-i} \geq \frac{r_k}{s-k}$. Suppose further that $r_k$ and $(s-k)$ are coprime. Then $f(x)$ has an irreducible factor of degree at least $s-k$ over $\Q$. 
\end{theorem}

In this paper, we show that the above coprimality condition of $r_k$ and $s-k$ can be weakened. Moreover, invoking the theory of Newton polygon (defined below), we improve the lower bound for the degree of an irreducible factor. Furthermore, we shall prove the theorem for polynomials having coefficients from the valuation ring of an arbitrary valued field and indicate how to derive it for polynomials with integer coefficients.\\

Let $v$ be a Krull valuation of a field $K$ with value group $G_v$ and valuation ring $R_v$ having maximal ideal $M_v$. We shall denote by $v^x$ the Gaussian prolongation of $v$ to $K(x)$ defined on $K[x]$ by 
\begin{equation}\label{eq:gauss}
v^x\big(\sum\limits_{i}a_ix^i\big) = \min\limits_{i}\{v(a_i)\}, ~~a_i \in K.
\end{equation}
For a polynomial $h(x)\in R_v[x]$, $\bar{h}(x)$ will stand for the polynomial over $R_v/M_v$ obtained by replacing each coefficient of $h(x)$ by its $v$-residue.  Below we define the notion of Newton polygon (see \cite[Section 6.4]{FQG}, \cite[Definition 1.D]{Jh-Kh} for details).
\begin{definition}\label{newton2}
Let $(K, v)$ be a valued field with value group $G_v$ and valuation ring $R_v$ having maximal ideal $M_v$. Let $\phi(x)\in R_v[x]$ be a monic polynomial with $\bar{\phi}(x)$ irreducible over $R_v/M_v$ and $v^x$ be the Gaussian prolongation defined by $(\ref{eq:gauss}).$ Let $f(x)\in R_v[x]$ be a polynomial not divisible by $\phi(x)$ with $\phi$-expansion\footnote{On dividing by successive powers of $\phi(x)$, every polynomial $f(x)\in K[x]$ can be uniquely written as a finite sum $\sum_{i\geq 0}f_i(x)\phi(x)^i$ with $\deg(f_{i}(x))<\deg(\phi(x))$, called the $\phi$-expansion of $f(x)$.} $\sum\limits_{i=0}^{n}a_i(x)\phi(x)^i, a_n(x) \neq 0$.  Let $P_i$ stand for the pair $(i, v^x(a_{n-i}(x)))$ when $a_{n-i}(x) \neq 0, ~0\leq i\leq n$. For distinct pairs $P_i, P_j,$ let $\mu_{ij}$ denote the element of the divisible closure of $G_v$ defined by \vspace*{-0.15in}$$\mu_{ij} = \dfrac{v^x(a_{n-j}(x)) - v^x(a_{n-i}(x))}{j-i};\vspace*{-0.15in}$$ in case $v$ is a real valuation, $\mu_{ij}$ is the slope of the line segment joining $P_i$ and $P_j$. Let $i_1$ denote the largest index $0 < i_1 \leq n$ such that 
\vspace*{-0.15in}
$$\mu_{0i_1} = \min\limits_{j} \left\lbrace \mu_{0j}: 0<j\leq n, a_{n-j}(x) \neq 0 \right\rbrace.\vspace*{-0.15in}$$
If $i_1 < n$, let $i_2$ be the largest index such that $i_1 < i_2 \leq n$ and \vspace*{-0.15in}$$ \mu_{i_1i_2}= \min\limits_{j}
  \left\lbrace \mu_{i_1 j}: i_1<j\leq n, a_{n-j}(x)\neq 0 \right\rbrace. \vspace*{-0.15in}$$
  Proceeding in this way if $i_{\y} = n$, then the $\phi$-Newton polygon of $f(x)$  is said to have $\y$ many edges whose slopes are defined to be $\mu_{0i_1}, \mu_{i_1i_2}, \cdots, \mu_{i_{\y-1}i_{\y}}$ which are in strictly increasing order. The pairs $P_0,P_{i_1},P_{i_2},\cdots, P_{i_{\y}}$ are called the successive vertices of the $\phi$-Newton polygon of $f(x)$ with respect to the valuation $v$.
\end{definition}
We now state our main result.
\begin{theorem}\label{thm2}
Let $v$ be a Krull valuation of a field $K$ with value group $G_v$ and valuation ring $R_v$ having maximal ideal $M_v$. Let $\phi(x)\in R_v[x]$ be a monic polynomial of degree $m$ which is irreducible modulo $M_v$. Let $f(x)\in R_v[x]$ be a polynomial not divisible by $\phi(x)$.
 Assume that the $\phi$-Newton polygon of $f(x)$ has $\y$ many edges  with positive slopes $\lambda_j$, $1\leq j\leq \y$. If $\x_j$ is the smallest positive number such that $\x_j\lambda_j \in G_v$ for $1\leq j\leq \y$, then $f(x)$ has an irreducible factor of degree at least $\max\limits_{1\leq j\leq \y}\{\x_jm\}$ over $K$.
\end{theorem}
With the notations as in the above theorem, the following corollary is an immediate consequence of Theorem \ref{thm2}.
\begin{corollary}\label{cor1.6}
Let the notations and assumptions be as in the above theorem. Then for any factorization $f_1(x)f_2(x)$ of $f(x)$ over $K$, we have
$$\min\{\deg f_1(x), \deg f_2(x)\} \leq \deg f(x) - \max\limits_{1\leq j\leq \y}\{\x_jm\}.$$
\end{corollary}

When $G_v = \Z$,  Theorem \ref{thm2} immediately gives the following result.  
\begin{theorem}\label{thm1}
Let $p$ be a prime number and $f(x) = \sum\limits_{i=0}^{n}a_ix^i$ with $a_0\neq 0$ be a polynomial having integer coefficients.  Assume that the $x$-Newton polygon of $f(x)$ with respect to the $p$-adic valuation $v_p$ has $\y$ many edges with positive slopes $\lambda_j$, $1\leq j\leq \y$. If $\lambda_j = \frac{r_j}{s_j}$ with $\gcd(r_j, s_j) = 1$, then $f(x)$ has an irreducible factor of degree at least $\max\limits_{1\leq j\leq \y}\{s_j\}$ over $\Q$.  
\end{theorem}

\begin{remark}
It may be noted that Theorem \ref{a1} follows from Theorem \ref{thm1} because in Theorem \ref{a1}, the $x$-Newton polygon of $f(x)$ with respect to $v_p$ has an edge having positive slope $\frac{r_k}{s-k}$ with $\gcd(r_k, s-k) = 1$.
\end{remark}

For the history of the problem and related literature, the reader may refer to \cite{Jak}. Before we get down to the proof of the theorem, an example may illustrate the importance of the result.
\begin{example}
Let $a, b$ be integers with $b \neq 0$. Let  $q$ be a prime number and $s, r$ denote the highest power of $q$ dividing $a, b$ respectively. Let $n, m$ be positive integers with $n > m$ and $(n-m)r > ns > 0$. Then  by Theorem \ref{thm1}, $f(x) = x^n + ax^m + b$ has an irreducible factor of degree at least $\max\{\frac{m}{\gcd(m, r-s)}, \frac{n-m}{\gcd{(n-m, s)}}\}$ over $\Q$.  In the particular cases, when either $n = m + 1$ with $\gcd(n-1, r-s) = 1$ or $m = 1$ with $\gcd{(n-1, s)} = 1$, then $f(x)$ has an irreducible of degree at least $n-1$ over $\Q$. Therefore, in these cases, if $f(x)$ does not have a linear factor, then $f(x)$ is irreducible over $\Q$. However, Theorem \ref{a1} provides no information about the irreducible factor of $f(x)$  when $n = m+1$.
\end{example}

\section{Proof of Theorem \ref{thm2}.}
\subsection{Some Notations and Definitions.}
Let $v$ be a Krull valuation of a field $K$ with value group $G_v$ and valuation ring $R_v$ having maximal ideal $M_v$.  We fix a prolongation $\tilde{v}$ of $v$ to an algebraic closure $\widetilde{K}$ of $K$. For an element $\alpha$ belonging to the valuation ring $R_{\tilde{v}}$ of $\tilde{v}$, $\bar{\alpha}$ will denote its $\tilde{v}$-residue, i.e., the image of $\alpha$ under the canonical homomorphism from $R_{\tilde{v}}$ onto its residue field $R_{\tilde{v}}/M_{\tilde{v}}.$

\begin{definition}\label{minimal}
Let $v$ be a henselian Krull valuation of a field $K$ and $\tilde{v}$   the unique prolongation of $v$ to the algebraic closure $\widetilde{K}$ of $K$ with value group $G_{\tilde{v}}$. A pair $(\alpha, \delta)$ belonging to $\widetilde{K}\times G_{\tilde{v}}$ is called a minimal pair (more precisely $(K, v)$-minimal pair) if whenever $\beta$ belongs to $\widetilde{K}$ with $[K(\beta) : K] < [K(\alpha) : K]$, then $\tilde{v}(\alpha - \beta) < \delta.$
\end{definition}
\begin{example}\label{kittu}
If $\phi(x)$ is a monic polynomial of degree $m\geq 1$ with coefficients in $R_v$ such that $\bar{\phi}(x)$ is irreducible over the residue field of $v$ and $\alpha$ is a root of $\phi(x)$, then $(\alpha, \delta)$ is a $(K, v)$-minimal pair for each positive $\delta$ in $G_{\tilde{v}}$, because whenever $\beta\in \widetilde{K}$ has degree less than $m$, then $\tilde{v}(\alpha - \beta) \leq 0$, for otherwise $\bar{\alpha} = \bar{\beta},$ which in view of the Fundamental Inequality \cite[Theorem 3.3.4]{En-Pr} would lead to $[K(\beta) : K] \geq m.$
\end{example}
\begin{definition}
Let $(K, v), (\widetilde{K}, \tilde{v})$ be as in the above definition and $(\alpha, \delta)$ belonging to $\widetilde{K}\times G_{\tilde{v}}$ be a $(K, v)$-minimal pair. The valuation $\widetilde{w}_{\alpha, \delta}$ of a simple transcendental extension $\widetilde{K}(x)$ of $\widetilde{K}$ defined on $\widetilde{K}[x]$ by
\begin{equation}\label{eq:a}
\widetilde{w}_{\alpha, \delta}\big(\sum\limits_{i}c_i(x-\alpha)^i\big) = \min\limits_{i}\{\tilde{v}(c_i) + i\delta\},~~ c_i \in \widetilde{K}
\end{equation}
will be referred to as the valuation with respect to the minimal pair $(\alpha, \delta);$ the restriction of $\widetilde{w}_{\alpha, \delta}$ to $K(x)$ will be denoted by $w_{\alpha, \delta}.$ For more details on the classification of extensions of a valuation from a base field $K$ to rational function fields, the reader may refer to the paper $\cite{FVK}$ and references therein. 
\end{definition}
With $(\alpha, \delta)$ as above, if $\phi(x)$ is the minimal polynomial of $\alpha$ over $K$, then it is well known \cite[Theorem 2.1]{A-P-Z} that for any polynomial $f(x)\in K[x]$ with $\phi$-expansion $\sum\limits_{i}a_i(x)\phi(x)^i$, one has
\begin{equation}\label{eq:b}
w_{\alpha, \delta}(f(x)) = \min\limits_{i}\{\tilde{v}(a_i(\alpha)) + iw_{\alpha, \delta}(\phi(x))\}.
\end{equation}
\begin{remark}
In particular,  if $(\alpha, \delta)$ is a minimal pair of the type described in Example \ref{kittu} with $\phi(x)$ as the minimal polynomial of $\alpha$ over $K$, then  for any polynomial $h(x) = \sum\limits_{i=0}^{m-1}a_ix^i\in K[x]$ having degree less than $m = \deg \phi(x)$, one has 
\begin{equation}\label{eq:c}
\tilde{v}(h(\alpha)) = v^x(h(x)).
\end{equation}
Clearly $(\ref{eq:c})$ needs to be verified when $m > 1$. Keeping in mind that $\bar{\phi}(x)$ is irreducible over $R_v/M_v$ of degree $m > 1$, it follows that $\tilde{v}(\alpha) = 0.$ If $(\ref{eq:c})$ were false, then the triangle inequality would imply that $\tilde{v}(h(\alpha)) > \min\limits_{i}\{\tilde{v}(a_i\alpha^i)\} = v(a_j)$ (say), which yields $\sum\limits_{i=0}^{m-1}{\big(\overline{\frac{a_i}{a_j}}\big)}(\bar{\alpha})^i = \bar{0}$, contradicting the fact that $\bar{\alpha}$ is a root of an irreducible polynomial $\bar{\phi}(x)$ of degree $m$. Hence $(\ref{eq:c})$ is proved.
\end{remark}
\begin{notation}\label{not}
Let $(\alpha, \delta)$ be as in Definition 2.3 and $\phi(x)$ be the minimal polynomial of $\alpha$ having degree $m$ over $K$. Let $w_{\alpha, \delta}$ be as in $(\ref{eq:b})$. For any non-zero polynomial $f(x)\in K[x]$ with $\phi$-expansion $\sum\limits_{i}a_i(x)\phi(x)^i$, we shall denote by $I_{\alpha, \delta}(f), S_{\alpha, \delta}(f)$ respectively the minimum and the maximum integers belonging to the set $\{i~|~w_{\alpha, \delta}(f(x)) = \tilde{v}(a_i(\alpha)) + iw_{\alpha, \delta}(\phi(x))\}$.
\end{notation}

With the above notation, the following result proved in \cite[Lemma 2.1]{Kh-Ku} will be used in the sequel.\\
\noindent \textbf{Theorem 2.A.}  For any non-zero polynomials $g(x), h(x)$ in $K[x]$, one has
 \begin{itemize}
 \item[(i)] $I_{\alpha, \delta}(g(x)h(x)) = I_{\alpha, \delta}(g(x)) + I_{\alpha, \delta}(h(x)),$
 \item[(ii)] $S_{\alpha, \delta}(g(x)h(x)) = S_{\alpha, \delta}(g(x)) + S_{\alpha, \delta}(h(x)).$
 \end{itemize}
 \subsection*{Proof of Theorem \ref{thm2}.}

Let $f(x) = \sum\limits_{i=0}^{n}a_i(x)\phi(x)^i$ be the $\phi$-expansion of $f(x)$. Let the set $\{(n-k_0, v^x(a_{k_0}(x))), (n-k_1, v^x(a_{k_1}(x))), \cdots, (n-k_{\y}, v^x(a_{k_{\y}}(x)))\}$ denote the successive vertices corresponding to all the edges having positive slopes in the $\phi$-Newton polygon of $f(x)$ where the $k_j$'s are integers with $k_{0}>k_1>\cdots>k_{\y}$. By the definition of the $\phi$-Newton polygon of $f(x)$, the last vertex of the Newton polygon is $(n, v^x(a_0(x)))$. So we have $k_l = 0$. Observe that the slope $\lambda_j$ is given by
\begin{equation}\label{2a}
\lambda_j =  \frac{v^x(a_{k_{j}}(x)) - v^x(a_{k_{j-1}}(x))}{k_{j-1} - k_{j}}, ~~1\leq j\leq \y.
\end{equation}
We may assume that $(K, v)$ is henselian because  the value group and the residue field remain the same on replacing $(K, v)$ by its henselization; moreover, if there is an irreducible factor of degree $\geq d$ in the factorization of $f(x)$ over the henselization, then there is an irreducible factor of degree $\geq d$ in the factorization of $f(x)$ over $K$. Let $\tilde{v}$ denote the unique prolongation of $v$ to the algebraic closure $\widetilde{K}$ of $K$.   
 Let $\alpha$ be a root of $\phi(x)$ in $\widetilde{K}$. Write $\phi(x) = c_m(x-\alpha)^m + \cdots + c_1(x-\alpha),~ c_m = 1.$ Define a positive element $\delta_j$ in the divisible closure $G_{\tilde{v}}$ of $G_v$ by 
$$\delta_j = \max\limits_{1\leq i\leq m}\bigg\{\frac{\lambda_j - \tilde{v}(c_i)}{i}\bigg\}.$$
Note that $\delta_j$ is positive in view of the fact that $c_m = 1$ and $\lambda_j > 0$. So $(\alpha, \delta_j)$ is a $(K, v)$-minimal pair in view of Example \ref{kittu}.
Let $\widetilde{w}_{\alpha, \delta_j}$ denote the valuation of $\widetilde{K}(x)$ defined by (\ref{eq:a}). Then by the choice of $\delta_j$, we have
$$\widetilde{w}_{\alpha, \delta_j}(\phi(x)) = \min\limits_{i}\{\tilde{v}(c_i) + i\delta_j\} = \lambda_j.$$ Keeping in mind $(\ref{eq:b})$ and $(\ref{eq:c})$, we see that
\begin{equation}\label{2b}
w_{\alpha, \delta_j}(f(x)) = \min\limits_{0\leq i\leq n}\{v^x(a_i(x)) + i\lambda_j\}.
\end{equation}
Let $I_{\alpha, \delta_j}(f)$ and $S_{\alpha, \delta_j}(f)$ be as in Notation \ref{not}. We claim that 
\begin{equation}\label{2c}
I_{\alpha, \delta_j}(f) = k_{j}, ~~~~~~~~S_{\alpha, \delta_j}(f) = k_{j-1}.
\end{equation}
Recall that $\lambda_j$ is the positive slope of the $\phi$-Newton polygon of $f(x)$ connecting the vertices $(n-k_{j-1}, v^x(a_{k_{j-1}}(x)))$, $(n-k_{j}, v^x(a_{k_{j}}(x)))$. By virtue of Definition \ref{newton2}, we see that

\begin{equation}\label{2e}
\min\limits_{0\leq i < k_{j-1}}\bigg\{\dfrac{v^x(a_i(x)) - v^x(a_{k_{j-1}}(x))}{k_{j-1} - i}\bigg\} \geq \dfrac{v^x(a_{k_{j}}(x)) -v^x(a_{k_{j-1}}(x)) }{k_{j-1} - k_{j}} = \lambda_j,
\end{equation} 
\begin{equation}\label{2f}
\min\limits_{k_{j-1} \leq i \leq n}\bigg\{\dfrac{v^x(a_{k_{j}}(x))-v^x(a_i(x))}{ i- k_{j}}\bigg\} \leq \dfrac{v^x(a_{k_{j}}(x)) -v^x(a_{k_{j-1}}(x)) }{k_{j-1} - k_{j}} = \lambda_j.
\end{equation}
Note that the smallest index $i$ for which equality in $(\ref{2e})$ holds is $i = k_{j}$. On the other hand, $i = k_{j-1}$ is the largest index such that equality holds in $(\ref{2f})$.  Therefore keeping in mind $(\ref{2b})$,  it follows that
\vspace*{-0.1in}\begin{equation}\label{2g}
w_{\alpha, \delta_j}(f(x)) = \min\limits_{0\leq i\leq n}\{v^x(a_i(x)) + i\lambda_j\} =  v^x(a_{k_j}(x)) + k_j\lambda_j = v^x(a_{k_{j-1}}(x)) + k_{j-1}\lambda_j \vspace*{-0.1in}
\end{equation} 
and $I_{\alpha, \delta_j}(f) = k_{j}, ~S_{\alpha, \delta_j}(f) = k_{j-1}$.

Let $f(x) = f_1(x)f_2(x)\cdots f_t(x)$ be the factorization of $f(x)$ into irreducible factors over $K$. Denote $I_{\alpha, \delta_j}(f_r)$ by $k_{j}^{(r)}$  and $S_{\alpha, \delta_j}(f_r)$ by $k_{j-1}^{(r)}$ for $1\leq r\leq t.$ Applying Theorem 2.A together with $(\ref{2c})$, we see that
$$k_{j} = k_{j}^{(1)} + \cdots + k_{j}^{(t)},~~~ k_{j-1} = k_{j-1}^{(1)} + \cdots + k_{j-1}^{(t)}.$$
Since $k_{j-1} > k_{j}$ and $k_{j-1}-k_{j} = k_{j-1}^{(1)} - k_{j}^{(1)} + \cdots + k_{j-1}^{(t)} - k_{j}^{(t)}$, it follows that  $k_{j-1}^{(r)} - k_{j}^{(r)} > 0$ for some $r$, $1\leq r\leq t$.  Without  loss of generality, we may assume that 
$$k_{j-1}^{(1)} - k_{j}^{(1)} > 0.$$
Let $f_1(x) = \sum\limits_{u=0}^{d_1}b_u(x)\phi(x)^u$ be the $\phi$-expansion of $f_1(x)$. Then we have
$$w_{\alpha, \delta_j}(f_1(x))  = v^x(b_{k_{j}^{(1)}}(x)) + k_{j}^{(1)}\lambda_j = v^x(b_{k_{j-1}^{(1)}}(x)) + k_{j-1}^{(1)}\lambda_j.$$
The above equality implies that $(k_{j-1}^{(1)} - k_{j}^{(1)})\lambda_j \in G_v$. Since $\x_j$ is the smallest positive element such that $\x_j\lambda_j \in G_v$, it follows that
$$(k_{j-1}^{(1)} - k_{j}^{(1)}) \geq \x_j.$$
As $S_{\alpha, \delta_j}(f_1) = k_{j-1}^{(1)}$, the above inequality shows that
$$\deg f_1(x) \geq k_{j-1}^{(1)}m \geq (k_{j-1}^{(1)} - k_{j}^{(1)})m \geq \x_j m.$$
As $j$ is arbitrary, we therefore conclude that $f(x)$ has an irreducible factor of degree at least $\max\limits_{1\leq j\leq \y}\{\x_jm\}$ over $K$. \\

\noindent\textbf{Acknowledgements.}  The authors are thankful to SERB MATRICS Project No. MTR /2017/001006 and IMSc for financial support. The authors express their gratitude to the anonymous referees for some important suggestions which improved the exposition of this paper.


 \end{document}